\documentclass[12pt,draftcls,onecolumn]{IEEEtran}

\pdfoutput=1
\usepackage{float}
\usepackage{amsfonts}
\usepackage{amsmath}
\usepackage{amsthm}
\usepackage{amssymb,latexsym}
\newtheorem{definition}{Definition} 
\newtheorem{lem}{Lemma}

\newtheorem{thm}{Theorem}
\newtheorem{prob}{Problem}
\usepackage{graphicx}
\usepackage{float}
\usepackage{subfig}
\usepackage{caption}
\usepackage{mathrsfs}
\usepackage{cite}
\begin{document}

\title{Computational Techniques for Reachability Analysis of Partially Observable Discrete Time Stochastic Hybrid Systems}
\author{Kendra Lesser, \IEEEmembership{Student Member, IEEE}, Meeko Oishi, \IEEEmembership{Member, IEEE}
\thanks{K. Lesser and M. Oishi are with the Department of Electrical and Computer Engineering, University of New Mexico, Albuqerque, NM 87131 USA; e-mail: \{lesser, oishi\}@unm.edu; Tel.:+1 505 353 2424; Fax: +1 505 277 0299.}
\thanks{This research was funded by National Science Foundation (NSF) Career Award CMMI-1254990, NSF Award CPS-1329878, NSA Science of Security Lablet at North Carolina State University (subaward to the University of New Mexico), and start-up funding from the Department of Electrical and Computer Engineering, University of New Mexico}}
\maketitle

\begin{abstract}

Reachability analysis of hybrid systems has  been used as a safety verification tool to assess offline whether the state of a system is capable of remaining within a designated safe region for a given time horizon.  Although it has been applied to stochastic hybrid systems, little work has been done on the equally important problem of reachability under  incomplete or noisy measurements of the state.  Further, there are currently no computational methods or results for reachability analysis of partially observable discrete time stochastic hybrid systems.  We provide the first numerical results for solving this problem, by drawing upon existing literature on continuous state  partially observable Markov decision processes (POMDPs).  We first prove that the value function for the reachability problem (with a multiplicative cost structure) is piecewise-linear and convex, just as for discrete state POMDPs with an additive cost function.  Because of these properties, we are able to extend existing point-based value iteration techniques to the reachability problem, demonstrating its applicability on a benchmark temperature regulation problem.
\end{abstract}


\section{Introduction}

Stochastic hybrid systems provide a modeling framework well-suited for a wide range of applications.  They allow for versatile dynamics that incorporate codependent discrete and continuous states, often exhibited in systems that may switch between different modes of operation, and account for probabilistic uncertainty in those dynamics.  Having such a flexible framework is particularly important in the context of safety verification, where the assessment of a system's ability to meet rigorous safety requirements must be as accurate as possible.  Indeed, reachability analysis (determining whether a system's state stays within a given safe region and/or reaches a desired target set within some finite time horizon) for hybrid systems has been studied extensively  \cite{oishi1}, \cite{Prand1}, \cite{Mitchell1},  \cite{Abate1}, \cite{summers}.  

Equally important to safety verification, however,  is the consideration of not only stochastic and complex dynamics, but also of noisy or incomplete measurements of the state.  While there has been some work on deterministic hybrid systems with incomplete information \cite{Verm2} or uncertain hybrid systems with the assumption of a worst-case disturbance \cite{Ghaemi2013}, reachability analysis of a partially observable stochastic hybrid system has been approached only recently \cite{Ding2013}, \cite{Lesser2013}, and only theoretically; there are currently no computational results for reachability analysis of \emph{partially} observable stochastic hybrid systems.

Computational results for reachability analysis of perfectly observable stochastic hybrid systems are also limited. The reachability problem for discrete time stochastic hybrid systems (DTSHS) is a multiplicative cost stochastic optimal control problem \cite{Abate1}, which can equivalently be formulated as a Markov decision process (MDP).  Solutions via dynamic programming produce a state-based feedback controller designed to optimize the system according to some cost function (see \cite{bertsekas2}).  Unfortunately, dynamic programming  requires evaluation of the value function over all possible states, which is infinite when those states are continuous.   Discretization procedures can be employed to impose a finite number of states, as in \cite{SoudjaniSiam13}, which presents a formal adaptive gridding procedure for verification of DTSHS.  Gridding methods are unfortunately subject to the ``curse of dimensionality'' and can lead to an unacceptable number of states that render the dynamic program impossible to implement.  Other approximate solution strategies include approximate dynamic programming, where the value function of the dynamic program is approximated by a set of basis functions, as in \cite{Kariot2013}.  Even so, current applications are limited to those with only a few discrete and continuous states.

The reachability problem for a partially observable DTSHS (PODTSHS) can similarly be formulated as a partially observable MDP (POMDP).  However, POMDPs are plagued by dimensionality on an even greater scale than MDPs.  The common approach to solving POMDPs is to replace the growing history of observations and actions by a sufficient statistic, often called the belief state, which, for a POMDP with an additive cost function, is the distribution of the current state conditioned on all past observations and actions  \cite{bertsekas2}.  This belief state is treated as the perfectly observed true state, and MDP solution methods can then be applied.  However, given a continuous state space, the belief state is now a continuous function defined over an infinite domain, and it is impossible to enumerate over all such functions.  Therefore the study of efficient, approximate solutions to POMDPs is essential.

Although finding the solution to a general POMDP is hard \cite{Lusena2001}, many algorithms for approximating solutions to finite state POMDPs have been developed.  These mainly rely on point-based value iteration (PBVI) schemes that only consider a subset of the belief space to update the value function (for a survey of PBVI algorithms, see \cite{shani13}).  Such methods must be tailored to continuous state POMDPs because of the dimensionality of the belief state.

Many existing methods for continuous state POMDPs assume the belief state is Gaussian, such as in  \cite{Brooks06}, \cite{Zhou2010}, and represent the belief state in a parameterized form which is then discretized and solved as a discrete state MDP.  For problems where the belief cannot be represented adequately as a single Gaussian, however, these technques are subject to the same curse of dimensionality as large discrete state MDPs.  Other methods use a Gaussian representation of the belief state to find \emph{locally} optimal solutions, either by parameterizing the value function \cite{Berg2012} or by assuming maximum-likelihood observations \cite{Platt2010} \cite{Erez2010}.  An extension of \cite{Platt2010} to non-Gaussian beliefs was presented in \cite{Platt2011}, where the belief states are estimated using sampling.  Another sampling-based method that allows for a non-Gaussian belief state is given by \cite{Thrun2000}, where the belief state is updated according to a particle filter, and Monte Carlo methods and nearest-neighbor approximations estimate the value function.

PBVI techniques have also been extended to the case of continuous states in \cite{Porta06}, which showed that for continuous states and discrete actions and observations, the value function remains piecewise-linear and convex (as was shown for discrete state POMDPs by \cite{sondik}).  These properties can be exploited to approximate the value function by a finite set of ``$\alpha$-functions,'' which are a function of the true state of the system, and represent the value of being in that state, including the future expected rewards assuming optimal actions are taken.  Further, by representing these $\alpha$-functions and the belief states as linear combinations of Gaussians, updating the belief state and value function can be done in closed form.  This technique was extended to hybrid domains, where the discrete mode is hidden and the belief state is a function only of the continuous variable \cite{Brunskill2010}.  The authors of \cite{Porta06} also showed that the belief state can be approximated using a particle filter rather than as a sum of Gaussians, and the continuous state PBVI method still applied.  

The reachability problem for PODTSHS further complicates the already difficult problem of solving continuous state POMDPs.  As was shown in both \cite{Ding2013} and \cite{Lesser2013}, the belief state of the PODTSHS is no longer just the conditional distribution of the current state of the system, but must also include the distribution of a binary variable indicating whether the state of the system has remained within a safe region up to the previous time step.  This, coupled with the stochastic hybrid system dynamics, makes representing the belief state as a single Gaussian impossible, and using sampling to update the belief can be expensive.  

Therefore, as the first investigation into approximate solutions to the reachability problem for PODTSHS, we consider continuous state PBVI techniques as in\cite{Porta06} and \cite{Brunskill2010}.  These techniques are amenable to stochastic hybrid dynamics, and have already been demonstrated as effective in hybrid domains with a hidden discrete state.  In this paper we present several contributions to the solution of safety verification problems for PODTSHS. First, we show that even with the multiplicative cost structure of the reachability problem, as in \cite{Abate1} and \cite{Lesser2013}, the value function is piecewise-linear and convex under the assumption of discrete actions and observations.  Further, the belief state, defined over a hybrid domain, and value function maintain the closedness property of the belief and value function updates, when they are represented as weighted sums of Gaussians.  Proving the preservation of these ``nice'' properties enables the application of existing POMDP solution techniques.  Second, we exploit the structure of the belief state and value function to extend the technique of \cite{Porta06} and \cite{Brunskill2010} to the reachability problem.  We outline a solution method, and demonstrate its effectiveness on a temperature regulation problem.

The rest of the paper is organized as follows.  Section \ref{PODTSHS}  defines a PODTSHS, and formulates the reachability problem.  Sections \ref{POMDP} and \ref{back3} provide an overview of POMDPs and their exact solution, and point-based value iteration techniques, respectively.  PODTSHSs and POMDPs are related in Section \ref{back4}.   Section \ref{PBVI} establishes properties of the value function, demonstrates how PBVI techniques can be used to solve the reachability problem for PODTSHS, and also provides a bound on the error introduced in approximating the true value function.  Section \ref{PBVI} also shows that the value function and belief updates preserve the Gaussian representation.  Section \ref{examples} provides numerical results using a benchmark temperature regulation problem, and discusses computational issues.  Section \ref{conc} provides concluding remarks and future directions.

\section{Background}\label{background}

\subsection{Reachability for PODTSHS}\label{PODTSHS}

A hybrid system is characterized by a set of both discrete and continuous states with interacting dynamics: the discrete state may affect the evolution of the continuous dynamics, and the continuous dynamics may affect when the discrete state changes.  In the case of a DTSHS, both the discrete and continuous dynamics may be characterized by stochastic kernels, the product of which determines the stochastic transition kernel governing the combined discrete/continuous state of the system.  We present a slightly modified definition of a DTSHS first introduced in \cite{Abate1}.

\begin{definition}\label{dtshs}(Discrete Time Stochastic Hybrid System $\mathcal{H}$). A DTSHS is a tuple $\mathcal{H} = (\mathcal{X},\mathcal{Q},\mathcal{U},T_x,T_q)$ where
	
\begin{enumerate}

\item
$\mathcal{X} \subseteq \mathbb{R}^n$ is a set of continuous states
\item
$\mathcal{Q} = \{q_1, q_2, ... q_{N_q}\}$ is a finite set of discrete states with cardinality $N_q$, with $\mathcal{S} = \mathcal{X} \times  \mathcal{Q}$ the hybrid state space
\item $\mathcal{U}$ is a compact Borel space which contains all possible control inputs affecting discrete and continuous state transitions
\item$T_x : \mathcal{B}(\mathbb{R}^n)\times\mathcal{Q}\times \mathcal{S} \times \mathcal{U} \rightarrow [0,1]$ is a Borel-measurable stochastic kernel which assigns a probability measure to $x_{t+1}$ 	
given $s_k=(x_t,q_t),u_t, q_{t+1}\,\forall\, t$: $T_x(dx_{t+1}\in B\mid q_{t+1}, s_t,u_t)$ where $B\in \mathcal{B}(\mathbb{R}^n)$, the Borel $\sigma$-algebra on $\mathbb{R}^n$
\item $T_q : \mathcal{Q}\times \mathcal{S} \times \mathcal{U} \rightarrow [0,1]$ is a discrete transition kernel assigning a probability distribution to $q_{t+1}$ given $x_t,q_t,u_t,\, \forall \, t$
\end{enumerate}
\end{definition} 

Kernels $T_x$ and $T_q$ can be combined for ease of notation to produce one hybrid state transition kernel, denoted $\tau(\cdot)$, given by:
\begin{equation}\label{statetrans}
\tau(ds' \mid s,u) = T_x(dx' \mid x, q, u, q')T_q(q'\mid x,q,u) \\
\end{equation}
The discrete state $q_{t+1}$ update depends on $q_t$, $x_t$ and $u_t$, and the continuous state $x_{t+1}$ update depends on $x_t$, $u_t$, and according to the specific problem may also be governed by $q_t$, $q_{t+1}$, or both.  For ease of notation we assume that the discrete state updates first, and the updated discrete state affects the continuous state, i.e. that $T_x(dx_{t+1} \mid x_t, u_t, q_{t+1})$, although modifying $T_x$ to include $q_t$ would not alter any subsequent results.  

For a PODTSHS, it is assumed that only an observation process is available to the controller, of the form $y_t = (y_t^x, y_t^q)$, where $y_t^x$ is associated with $x_t$, and $y_t^q$ with $q_t$.  While $y_t^x$ could be continuous, for computational purposes we assume that it is discrete-valued, even though $x_t$ is continuous (which could arise simply by discretizing the observation process).  The observation process is given by
\begin{align}
y_t^x &= h(x_t,u_{t-1}) + v_t \label{obs1}\\
y_t^q &\sim Q_{q,y^q}(u) \label{obs2}
\end{align}
The probability that $y_t^q = n$, $P[y_t^q = n \mid q_t=q,\,u_{t-1}=u] = Q_{q,n}(u)$, is given by the state transition matrix $Q(u)$ which is dependent on the control input $u$.  For the continuous state observation $y_t^x$ that is continuous-valued, it is subject to additive noise $v_t$, which is independent and identically distributed with positive density $\varphi(v)$ (i.e. Gaussian), and the function $h$ is assumed to be bounded and continuous.  Otherwise we assume $y_t^x$ has a state transition matrix similar to $Q$, and we will write $\varphi(y^x\mid x,u)$ to express the conditional discrete distribution of $y^x$.  The filtrations $\mathcal{G}_t$ and $\mathcal{Y}_t$ are generated by the sequences $\{s_0,\dotsc,s_{t},y_1,\dotsc,y_{t-1}\}$ and $\{y_1,\dotsc,y_t\}$, respectively.  We also assume an initial Borel-measurable density on $s_0=(x_0,q_0)$, $s_0\sim\rho(x,q)\in P(\mathcal{S})$, i.e. that $\rho$ lies in the space of all probability measures on $\mathcal{S}$.  Finally, based on $\rho$, $\tau$, $\varphi$, and $Q(u)$, the probability measure $\mathbb{P}^{\pi}$ is induced by the control policy $\pi$ defined over the full state space $\Omega$, which includes $s_t$ and $y_t$ for all $t$.  

Next, we present a cost function to analyze the reachability of the partially observable DTSHS, i.e. the ability of the state to remain within some safe or desired region of the state space.  We want to find both a control policy that maximizes the probability of the state remaining within that desired set, as well as an estimate of that probability.  As in \cite{Abate1}, this problem can be formulated as a stochastic optimal control problem.  For a Borel set $K\subseteq \mathcal{X}\times \mathcal{Q}$, terminal time $T$, and predefined policy $\pi$, define the cost function as
\begin{equation}\label{probRA}
r_{K}(\pi) = \mathbb{P}^{\pi}[s_t \in K \, \forall \,t=0,\dotsc,T]
\end{equation}
Since for a random variable $X$, $\mathbb{P}[x \in A] = \mathbb{E}[{\bf 1}_A(x)]$, with $\mathbb{E}$ denoting expected value and indicator function ${\bf1}_A(x)=1$ if $x \in A$ and ${\bf1}_A(x)=0$ otherwise, \eqref{probRA} is rewritten as in \cite{Abate1}:
\begin{equation} \label{ERA}
r_K(\pi) ={\mathbb E}^{\pi} \left[ \prod_{t=0}^{T}{\bf 1}_K(s_t) \right]
\end{equation}
The expected value is taken with respect to the measure $\mathbb{P}^{\pi}$, hence the notation $\mathbb{E}^{\pi}$.  We want to maximize $r_K(\pi)$ with respect to the control policy $\pi$.  The set $\Pi$ of admissible policies will be restricted to non-randomized policies, i.e. in which $\pi(y_t)$ generates one control input $u_t$ with probability $1$.   The optimal policy $\pi^*$ is then given by
\begin{equation}\label{optpol}
\pi^* = \arg\sup_{\pi\in\Pi} \left\{r_K(\pi)\right\}
\end{equation}
We can now formally define the problem we wish to solve.

\begin{prob}\label{probState}
Consider a DTSHS $\mathcal{H}$  (defined in Definition \ref{dtshs}) with observations \eqref{obs1} - \eqref{obs2} and initial distribution $\rho(x,q)\in P(\mathcal{S})$.  Given a safe set $K$ and time horizon $T$ we would like to
\begin{enumerate}
\item
Compute the maximal probability of remaining within $K$ for $T$ time steps, given by $\sup_{\pi} r_K(\pi)$.
\item
Compute the optimal policy $\pi^*$ such that  $\sup_{\pi} r_K(\pi) = r_K(\pi^*)$.
\end{enumerate}
If the maximal probability and optimal policy cannot be computed exactly (which is quite likely \cite{Lusena2001}), an approximation producing a suboptimal policy and lower bound on the maximal reachability probability are desired.
\end{prob}

\subsection{Optimal Control of POMDPs}\label{POMDP}

POMDPs provide a framework for analyzing a discrete time system whose state depends on the actions of an agent (controller), who is trying to drive the state to optimize some objective. The state evolves stochastically and is Markovian (the state at the next time step depends only on the current state and action).  Further, in choosing actions, the agent can not directly observe the state of the system, instead only having access to an observation process.   We first define a POMDP with discrete states, actions, and observations, and an additive cost function.  The theory and solution techniques for this type of POMDP provide the foundation for our extension to a PODTSHS and the solution of Problem \ref{probState}.  

\begin{definition}\label{pomdpDef}(POMDP $\mathcal{G}$) A POMDP is a tuple $\mathcal{G} = (\mathcal{S}, \mathcal{U}, \mathcal{Y}, \tau, \psi, R)$ where

\begin{enumerate}
\item
$\mathcal{S}$ is a set of discrete states
\item
$\mathcal{U}$ is a discrete set of possible actions the agent can take
\item
$\mathcal{Y}$ is a set of discrete observations
\item
$\tau: \mathcal{S}\times\mathcal{S}\times\mathcal{U}\rightarrow [0,1]$ is a state transtion function assigning a probability distribution to state $s_{t+1}$ given state $s_t$ and  action $u_t$ for all $t$, $\tau(s_{t+1}\mid s_t, u_t)$
\item
$\psi: \mathcal{Y}\times\mathcal{S}\times\mathcal{U}\rightarrow [0,1]$ is an observation function assigning a probability distribution to observation $y_t$ given state $s_t$ and action $u_t$ for all $t$, $\psi(y_t\mid s_t, u_t)$
\item
$R: \mathcal{S}\times \mathcal{U}\rightarrow \mathbb{R}$ is a function assigning a reward (which we define as being in the set of all real numbers $\mathbb{R}$, although this could be generalized to any space)  at each time step $t$, given the current state $s_t$ and action $u_t$, $R(s_t, u_t)$
\end{enumerate}

\end{definition}

The goal for the POMDP $\mathcal{G}$ is to maximize the expected sum of rewards over a (possibly infinite) time horizon $T$ by optimally choosing a sequence of control actions $\overline{u}=\{u_1, u_2,\dotsc\}$.
\begin{equation}\label{sumCost}
\max_{\overline{u}} \mathbb{E}\left[\sum_{t=0}^T R(s_t, u_t)\right]
\end{equation}

Rather than keeping track of all past observations and actions in order to make an optimal decision at time $t$, a belief state is used instead, which summarizes all available information up to time $t$.  The belief state is a \emph{sufficient statistic} for the set of all observations and actions $\{u_1,\dotsc, u_{t-1}, y_1,\dotsc,y_t\}$ because it condenses all information necessary for making optimal decisions \cite{bertsekas2}.  In the case of an additive cost POMDP, the belief state is a probability density function that describes the probability of being in state $s$ given all past observations and actions, $b(s_t) = P[s_t\mid u_1,\dotsc, u_{t-1}, y_1,\dotsc,y_t]$.  Treating the belief state as the true state of the system, $\mathcal{G}$ can be equivalently solved as a perfect state information MDP.  An optimal policy $\pi^*$ for the POMDP is defined in terms of the belief state, and maps beliefs to actions: $\pi^*:\mathcal{B}\rightarrow \mathcal{U}$.

The optimal policy can be found by using a value function over the space of beliefs $\mathcal{B}$, which describes the cumulative reward from time $t$ to the final time $T$ (or over $T-t+1$ time steps), for a particular belief state $b$, and assuming the system behaves optimally from time $t+1$ to $T$.  The control $u_t$ is chosen to maximize the value function at a specific belief $b$.  Because the value function assumes only optimal actions are taken starting at time $t+1$, it can be defined recursively using the optimal value function at time $t+1$.
\begin{equation}
V_t^*(b) = \max_{u\in\mathcal{U}}\left \{\sum_s R(s,u)b(s) + \sum_y V_{t+1}^*(M_{y,u}[b])P[y\mid u, b]\right\}
\end{equation}
The transition operator $M_{y,u}[b]$ provides the next belief state $b_{t+1}$ given the current observation, action, and belief state.  Sondik \cite{sondik} first showed that for a finite horizon $T < \infty$, the value functions are piecewise-linear and convex, and thus can be expressed as
\begin{equation}\label{Valpha}
V_t^*(b) =\max_{\alpha_t^i\in\Gamma_t}\, \sum_s \alpha_t^i(s) b(s)
\end{equation}

The functions $\alpha_t^i\in\Gamma_t$, or ``$\alpha$-vectors'', can be thought of as representing a policy tree starting from a specific action $u$ and state $s$, which then specifies optimal actions conditioned on observations for the following time steps $t+1$ to $T$.  The $\alpha$-vectors thus characterize the current value of being in state $s$ and taking action $u$, plus the expected sum of future rewards assuming all subsequent actions are chosen optimally.  Because each $\alpha$-vector is associated with a specific action, by picking the $\alpha$-vector that maximizes $\sum_s \alpha_t^i(s) b(s)$, we are also defining the optimal policy for belief $b$ at time $t$.

In order to calculate the value function and optimal policy for all times $t$, all that is required are the complete sets of $\alpha$-vectors, $\Gamma_t$, for all $t$.  Unfortunately, the number of $\alpha$-vectors grows exponentially with $t$.  The $\alpha$-vectors at time $t$ are computed recursively from the $\alpha$-vectors calculated at time $t+1$.  For each action, we observe one of $|\mathcal{Y}|$ observations (where $|\cdot|$ indicates the cardinality of the set), and for each of those observations there is a subsequent $\alpha$-vector defined at time $t+1$, resulting in $|\mathcal{U}||\Gamma_{t+1}|^{|\mathcal{Y}|}$ $\alpha$-vectors at time $t$.  

Often, some of the $\alpha$-vectors are completely \emph{dominated} by another $\alpha$-vector or set of $\alpha$-vectors (where $\sum_s \alpha_t^j(s)b(s) < \sum_s \alpha_t^k(s)b(s)$ for all $b(s)$ implies $\alpha_t^j$ is dominated by $\alpha_t^k$).  While those dominated vectors clearly do not need to be included in the set $\Gamma_t$, finding the unnecessary $\alpha$-vectors is also computationally expensive.  A number of approximate solution techniques, including point-based value iteration, have been developed.

\subsection{Point-Based Value Iteration}\label{back3}

Point-Based Value Iteration (PBVI) computes the value function only over a finite subset $B\subset \mathcal{B}$.  The  general idea is to generate a collection of points $b \in \mathcal{B}$, and for each of these points perform a ``backup'' operation to get a new estimate of the value function at that point.  Most PBVI approaches use the same method of updating the value function at each belief point (the ``backup'' operation) and are distingushed by how they select the subset $B$ (see \cite{shani13}).  Here we outline the method of estimating the value function presuming a set $B$ has already been selected.  A discussion of various methods for selecting $B$ can be found in \cite{pineau06} and \cite{shani13}.

One $\alpha$-vector must be generated for each belief point $b^i\in B$,  $B = (b^0, b^1,\dotsc, b^n)$, so that $\tilde{\Gamma}_t = (\alpha_t^0, \alpha_t^1, \dotsc, \alpha_t^n)$ for all $t$.  We assume that an $\alpha$-vector $\alpha_t^j$ corresponding to $b^j$  will apply to all belief points in a region around $b^j$ (i.e. for any $b$ in a neighborhood of $b^j$ the same action will likely be optimal).  Hence the value at some $b$ not necessarily in $B$ can be approximated by
\begin{equation*}
V_t^*(b) \approx \max_{\alpha_t^i\in\tilde{\Gamma}_t}\, \sum_s \alpha_t^i(s) b(s)
\end{equation*}
as in \eqref{Valpha} but with a restricted set $\tilde{\Gamma}_t\subset \Gamma_t$.  The set $\tilde{\Gamma}_t$ is generated recursively from $\tilde{\Gamma}_{t+1}$, but without enumeration over all possible combinations of observations and subsequent $\alpha$-vectors in $\tilde{\Gamma}_{t+1}$ (the full policy tree starting at time $t$).  

For a specific $b \in B$, the value function at time $t$ can be approximated as follows (see, e.g., \cite{shani13} for more detail):
\begin{align}
V_t^*(b) &= \max_{u\in\mathcal{U}} \left\{ \sum_{s\in\mathcal{S}} R(s,u)b(s) + \sum_y V_{t+1}^*(M_{y,u}[b])P[y\mid u, b]\right\} \\
&=  \max_{u\in\mathcal{U}} \left\{ \sum_{s\in\mathcal{S}} R(s,u)b(s) + \sum_y \max_{\alpha_t^i\in\tilde{\Gamma}_t} \sum_{s'\in\mathcal{S}} \alpha_{t+1}^i(s')M_{y,u}[b](s') P[y\mid u,b] \right\}\\
&=  \max_{u\in\mathcal{U}} \left\{ \sum_{s\in\mathcal{S}} R(s,u)b(s) + \sum_y \max_{\alpha_{t+1}^i\in\tilde{\Gamma}_{t+1}} \sum_{s'\in\mathcal{S}}\alpha_{t+1}^i(s')\psi(y\mid s',u)\sum_s\tau(s'\mid s,u) b(s) \right\} \label{opexpand}
\end{align}
where \eqref{opexpand} follows from expanding the operator $M$.  We define the restricted set $\tilde{\Gamma}_t$ using the following expressions
\begin{align}
\alpha_{y,u}^i(s) &= \sum_{s'\in\mathcal{S}}\alpha_{t+1}^i(s')\psi(y\mid s',u)\tau(s'\mid s,u) \label{alpha1} \\
\alpha_{y,u,b}(s) &= \arg\max_i\sum_s \alpha_{y,u}^i(s)b(s) \label{alpha2}
\end{align}
to obtain 
\begin{equation}\label{alpha_set}
\tilde{\Gamma}_t = \bigcup_{b\in B}\left\{R(s,u) + \sum_y\alpha_{y,u,b}(s) \right\}_{\forall u \in \mathcal{U}}
\end{equation}
The function \eqref{alpha1} is an $\alpha$-function corresponding to a specific action $u$ and observation $y$ (representing the value of being in state $s$ given $y$ is observed and action $u$ is taken).  For a given belief state $b$, \eqref{alpha2} is the optimal function $\alpha_{y,u}^i$ for that belief state given $y$ is observed and action $u$ is taken.  Summing over all observations $\mathcal{Y}$ (essentially taking the expected value with respect to $y$) and then taking the union over all belief states in $B$ and actions $u\in\mathcal{U}$ produces \eqref{alpha_set}, the set of $\alpha$-functions at time $t$.

Finally, we define the backup operator for a specific belief point $b$ using the set $\tilde{\Gamma}_t$ \eqref{alpha_set} as
\begin{equation}\label{backup}
backup(b) = \arg\max_{\alpha_t^i\in\tilde{\Gamma}_t} \sum_{s\in\mathcal{S}} \alpha_t^i(s)b(s)
\end{equation}
Note that we can now define the optimal value function \eqref{opexpand} as 
\begin{equation}\label{backupVal}
V_t^*(b) = \sum_s( backup(b)(s)\times b(s))
\end{equation}
The overall PBVI algorithm then consists of selecting a set of belief points $B$, and repeatedly applying \eqref{backup} to each element of $B$.  In the case of a finite horizon of length $T$, the backup operator will be applied $T$ times, and for an infiite horizon, the backup operator will be applied until some tolerance level is reached (for example, where $\|V_{n+1}(b) - V_n(b)\| < \epsilon$).

The above derivations apply to a model with discrete state, action, and observation spaces, but \cite{Porta06} actually shows that the same technique applies to a POMDP with a continuous state space and discrete observation and action spaces.  In this case, the $\alpha$-vectors are replaced by $\alpha$-functions defined over the continuous space $S$.  Because the observations and actions are assumed discrete, there are a finite number of these $\alpha$-functions, and so the value function is still piecewise-linear and convex, but now with respect to the $\alpha$-functions.  In this case, the optimal value function may instead be represented as $V_t^*(s) = \sup_{\alpha_t^i\in\tilde{\Gamma}_t}\int_{\mathcal{S}} \alpha_t^i(s) b(s)\,ds$.

 When replacing $\mathcal{S}$ with a continuous state space, all of the above derivations hold, but all summations over $\mathcal{S}$ are replaced by integrals.  To generalize from the purely discrete case, \cite{Porta06} uses inner product notation rather than a summation or integral, so that \eqref{backup} would instead be written as
\begin{equation*}
backup(b) = \arg\max_{\alpha_t^i\in\tilde{\Gamma}_t}\, \langle \alpha_t^i,b\rangle
\end{equation*}
We maintain this notation in our derivations, where in the case of a hybrid state space with continuous state $x$ and discrete state $q$, $\langle f,g\rangle = \sum_q \int f(x,q)g(x,q)\,dx$ for well-defined functions $f$ and $g$.

 All of our derivations will assume discrete actions, and discrete (or discretized) observations.  If the assumption of discrete actions and observations is dropped, the value function is still convex but is no longer piecewise-linear (since there are an infinite number of $\alpha$-functions at any given time step).  The authors of \cite{Porta06} show, however, that a PBVI algorithm can still be applied to estimate the value functions by carefully sampling from the observation and action spaces.  While our method can also be extended to continuous actions and observations, we assume they are discrete for clarity and completeness of subsequent derivations.

\subsection{Relating Problem \ref{probState} to a POMDP}\label{back4}

We write the PODTSHS of Problem \ref{probState} as a  POMDP, which we denote $\mathcal{G}-hybrid$, with hybrid state space $\mathcal{S} = \mathcal{X}\times\mathcal{Q}$, control space $\mathcal{U}$, hybrid observation space $\mathcal{Y} = \mathcal{Y}^x\times\mathcal{Y}^q$, state transition function $\tau$ given by \eqref{statetrans}, and observation model $\psi(y\mid s,u) = Q_{q,y^q}(u)\varphi(y^x - h(x,u))$.  The reward function is given by $R(s_t,u_t) = {\bf 1}_K(s_t)$.  Note, however, that in contrast to the maximization over a sum of $R(s_t, u_t)$ as in \eqref{sumCost} for POMDP $\mathcal{G}$, we want to maximize the product for $\mathcal{G}-hybrid$, as described in Problem \ref{probState}.

We then reformulate $\mathcal{G}-hybrid$ into an equivalent perfect state information MDP, in the same fashion as for POMDP $\mathcal{G}$, by redefining the state of the system in terms of a sufficient statistic, or belief state.  However, because the cost function \eqref{ERA} is multiplicative rather than additive, the posterior distribution of the state at time $t$ given all available information up to time $t$ is no longer valid.  In \cite{Lesser2013}, we developed an appropriate sufficient statistic to solve \eqref{ERA} as a perfect state information problem using standard dynamic programming techniques. 

In summary, a change of measure, $\mathbb{P}^{\dag}$, makes the observation processes $\{y_t^x\}$ and $\{y_t^q\}$ each identically distributed and independent of $\{x_t\}$ and $\{q_t\}$, respectively, via the Radon-Nikodym derivative \cite{Lesser2013} \cite{SteinShak}, such that
\begin{equation}\label{RNder}
\left. \frac{d\mathbb{P}^{\pi}}{d\mathbb{P}^{\dag}}  \right|_{\mathcal{G}_t} = \Lambda_t
\end{equation}
where
\begin{equation*}
\Lambda_t =  \prod_{l=1}^t \frac{\varphi(y_l^x - h(x_l,u_{l-1}))Q_{q_l,y_l^q}(u_{l-1})}{\varphi(y_l^x)\frac{1}{N_q}}
\end{equation*}

The change of measure facilitates sampling to generate the belief states.  The sufficient statistic, $\sigma(x,q)$, can be defined as 
\begin{equation}\label{sigma}
\sigma_t(x, q) = \mathbb{E}^{\dag}\left[\left. {\bf 1}_q(q_t){\bf 1}_{x}(x_t) \prod_{i=0}^{t-1}{\bf 1}_K(s_i) \Lambda_t \right| \mathcal{Y}_t\right],
\end{equation}
a modification of the posterior distribution, that represents an unnormalized conditional density of the current state joined with the probability that all previous states are in $K$.  The sufficient statistic can be updated recursively using a bounded linear operator $\Phi$:
\begin{equation}\label{sigma_rec}
\begin{cases}
\sigma_0(x,q) = \rho(x,q) \\
\sigma_t(x,q) = \Phi_{y,u}[\sigma_{t-1}](x,q)
\end{cases}
\end{equation}
where $\Phi_{y,u}[\sigma]$ is given by
 \begin{equation}\label{sigma_op}
\Phi_{y,u}[\sigma](x',q') = \sum_{q\in\mathcal{Q}}N_{y^q}N_{y^x}Q_{q',y^q}(u)\int_{\mathbb{R}^n} {\bf 1}_K(x,q)  \varphi(y^x\mid x',u) \tau(x',q'\mid x,q,u)\sigma(x,q)\,dx
\end{equation}
 in the case of discrete observations $y^x$, with $N_{y^q}$ the number of possible observations of discrete mode $q$, and $N_{y^x}$ the number of possible observations of continuous state $x$.

The dynamic programming recursion to solve for \eqref{ERA} and \eqref{optpol}, 
\begin{equation}\label{optDP}
\begin{cases}
V_T^*(\sigma) = \langle \sigma, {\bf 1}_K \rangle \\
V_t^*(\sigma) = \sup_{u\in\mathcal{U}} \mathbb{E}^{\dag}\left[ V_{t+1}^*(\Phi_{y,u}[\sigma])\right]
\end{cases}
\end{equation}
first evaluates the value function $V_T^*(\sigma)$ in terms of the sufficient statistic $\sigma$, then recursively solves $V_{T-1}^*(\sigma)$, $V_{T-2}^*(\sigma)$, etc., ultimately resulting in $V_0^*(\rho) = \sup_{{\pi}\in\Pi}r_K(\pi)$ (see \cite{Lesser2013} for proof that this is true).  We note that \cite{Ding2013}, \cite{Tkachev2013} showed that the reachability problem can be equivalently formulated as an additive cost optimization by modifying the state of the system to include a binary variable indicating whether the state has remained within the safe region up to the previous time.  The authors of \cite{Ding2013} developed and then used this additive cost formulation to generate the sufficient statistic for a partially observable DTSHS as the posterior distribution of the modified state.  In \cite{Lesser2013}, we showed its equivalence to the multiplicative cost formulation and sufficient statistic.

We write the recursive relationship between the value functions  using operator notation,
\begin{equation}\label{DPop}
V^*_{t} = H[V_{t+1}^*]
\end{equation}
with $H[V] = \sup_{u\in\mathcal{U}} \mathbb{E}^{\dag}\left[ V(\Phi_{y,u}[\sigma])\right]$ as in \eqref{optDP}.  A useful property of $H$, which we will use later, is that it is a nonexpansion, meaning

\begin{equation}\label{nonExp}
\|H[V] - H[U]\|_{\infty} \leq \|V - U\|_{\infty}
\end{equation}
The proof of \eqref{nonExp} is straightforward, and hence omitted.

Similarly to the POMDP $\mathcal{G}$, the value function in \eqref{optDP} for $\mathcal{G}-hybrid$ must be solved for all functions $\sigma$, which lie in an infinite dimensional space. This clearly cannot be solved directly.  However, we will show that $\mathcal{G}-hybrid$ maintains the properties of the POMDP $\mathcal{G}$, i.e. that the value function is piecewise-linear and convex, and can be expressed as in \eqref{Valpha}, but with a hybrid state $s$.  In turn, we can use PBVI techniques to approximate the solution to Problem \ref{probState}.

\section{Point-Based Value Iteration for Hybrid Dynamics and Multiplicative Cost}\label{PBVI}
\subsection{Properties of the Value Function}
We first demonstrate that the value function for Problem \ref{probState} is convex for a hybrid state space with possibly continuous (or hybrid) actions and observations, and that the value function is also piecewise-linear in the case of purely discrete actions and observations. 

\begin{lem}\label{lemma1}
The value function \eqref{optDP} is convex in $\sigma$ for all $k$.
\end{lem}
\begin{IEEEproof}
By induction, at time $T$ for $0\leq\lambda\leq1$
\begin{align*}
V_T^*(\lambda\sigma_1 + (1-\lambda)\sigma_2) &= \sum_{q\in\mathcal{Q}} \int_{\mathbb{R}^n} {\bf 1}_K(x,q) \left[\lambda\sigma_1(x,q) + (1-\lambda)\sigma_2(x,q)\right]\,dx \\
&= \lambda V_T^*(\sigma_1) + (1-\lambda)V_T^*(\sigma_2)
\end{align*}
Assuming $V_{t+1}^*(\sigma)$ is convex in $\sigma$ 
\begin{align*}
V_t^*(\lambda\sigma_1 + (1-\lambda)\sigma_2) &= \sup_{u\in\mathcal{U}} \sum_{y^q}\int_{\mathbb{R}^n} V_{t+1}^*(\Phi_{y,u}[\lambda\sigma_1 + (1-\lambda)\sigma_2])\frac{1}{N_q}\varphi(y^x)\,dy^x \\
&=\sup_{u\in\mathcal{U}} \sum_{y^q}\int_{\mathbb{R}^n} V_{t+1}^*(\lambda \Phi_{y,u}[\sigma_1] + (1-\lambda) \Phi_{y,u}[\sigma_2]) \frac{1}{N_q}\varphi(y^x)\,dy^x \\
&\leq \sup_{u\in\mathcal{U}} \sum_{y^q}\int_{\mathbb{R}^n}  \left[\lambda V_{t+1}^*(\Phi_{y,u}[\sigma_1]) + (1-\lambda) V_{t+1}^*(\Phi_{y,u}[\sigma_2])\right] \frac{1}{N_q}\varphi(y^x)\,dy^x \\
&\leq  \sup_{u\in\mathcal{U}} \sum_{y^q}\int_{\mathbb{R}^n}\lambda V_{t+1}^*(\Phi_{y,u}[\sigma_1])  \frac{1}{N_q}\varphi(y^x)\,dy^x\\
\end{align*}
\begin{align*}
&\hspace{10 mm} +  \sup_{u\in\mathcal{U}} \sum_{y^q}\int_{\mathbb{R}^n}(1-\lambda) V_{t+1}^*(\Phi_{y,u}[\sigma_2])   \frac{1}{N_q}\varphi(y^x)\,dy^x \\
&\leq \lambda V_t^*(\sigma_1) + (1-\lambda)V_t^*(\sigma_2)
\end{align*}
\end{IEEEproof}

\begin{lem}\label{lemma2}
For any $t$, the value function \eqref{optDP} can be written as
\begin{equation*}
V_t^*(\sigma) = \sup_{\alpha_t^i\in\Gamma_t}\, \langle \alpha_t^i, \sigma \rangle
\end{equation*}
\end{lem}
\begin{IEEEproof}
By induction, at time $T$
\begin{equation*}
V_T^*(\sigma) = \sum_{q\in\mathcal{Q}} \int_{\mathbb{R}^n}  {\bf 1}_K(x,q) \sigma(x,q) \,dx
\end{equation*}
By defining $\alpha_T(x,q) = {\bf 1}_K(x,q)$, 
we obtain the desired result.  Note that this definition of $\alpha_T$ is in line with the definition given in Section \ref{back3}, because although it does not represent a full policy tree (being at the terminal time, there are no more branches on the tree), it does represent the immediate value of being in state $(x,q)$, given by ${\bf 1}_K(x,q)$.

 Next, assuming $V_{t+1}^*(\sigma) = \sup_{\Gamma_{t+1}} \,\langle \alpha_{t+1}^i, \sigma \rangle$, $V_t^*$ can be written as

\begin{align*}
V_t^*(\sigma) &= \sup_{u\in\mathcal{U}} \sum_{y^q} \int_{\mathbb{R}^n}V_{t+1}^*(\Phi_{y,u}[\sigma]) \frac{1}{N_q} \varphi(y^x)\,dy^x \\
&=  \sup_{u\in\mathcal{U}} \sum_{y^q} \int_{\mathbb{R}^n} \sup_{\Gamma_{t+1}} \langle \alpha_{t+1}^i, \Phi_{y,u}[\sigma]\rangle  \frac{1}{N_q} \varphi(y^x)\,dy^x \\
&=  \sup_{u\in\mathcal{U}} \sum_{y^q} \int_{\mathbb{R}^n} \sup_{\Gamma_{t+1}}  \sum_{q'}\int_{\mathbb{R}^n} \alpha_{t+1}^i(x',q') \Phi_{y,u}[\sigma](x',q')\,dx' \frac{1}{N_q} \varphi(y^x)\,dy^x \\
&= \sup_{u\in\mathcal{U}} \sum_{y^q} \int_{\mathbb{R}^n} \sup_{\Gamma_{t+1}}  \sum_{q'}\int_{\mathbb{R}^n} \sum_q \int_{\mathbb{R}^n}  \alpha_{t+1}^i(x',q') Q_{q',y^q}(u) \varphi(y^x - h( x',u)) {\bf 1}_K(x,q) \\
&\hspace{60 mm} \times  \tau(x',q'\mid x,q,u)\sigma(x,q)\,dx\,dx'\,dy^x \\
&=  \sup_{u\in\mathcal{U}} \sum_{y^q} \int_{\mathbb{R}^n} \sup_{\Gamma_{t+1}} \sum_q\int_{\mathbb{R}^n}\left[ \sum_{q'}\int_{\mathbb{R}^n} \alpha_{t+1}^i(x',q') Q_{q',y^q}(u) \varphi(y^x - h( x',u)) \right.  \\
&\hspace{60 mm}\times \left. \vphantom{ \sum_{q'}\int_{\mathbb{R}^n}} \tau(x',q'\mid x,q,u)\,dx' \right] {\bf 1}_K(x,q) \sigma(x,q) \,dx\,dy^x \\
\end{align*}
\begin{align*}
&= \sup_{u\in\mathcal{U}} \sum_{y^q} \int_{\mathbb{R}^n} \sup_{\Gamma_{t+1}}\,\left \langle  \sum_{q'}\int_{\mathbb{R}^n} \alpha_{t+1}^i(x',q') Q_{q',y^q}(u) \varphi(y^x - h( x',u))\right. \\ 
&\hspace{60mm}\left. \vphantom{\int_{\mathbb{R}^n}} \times \tau(x',q'\mid x,q,u)\,dx' {\bf 1}_K(x,q), \sigma (x,q)\right \rangle \,dy^x
\end{align*}
Then for a specific observation $y$, action $u$, and $\alpha_{t+1}^i$ function, the function $\alpha_{y,u}^i$ can be defined as
\begin{equation}\label{alpha_yu}
\alpha_{y,u}^i(x,q) =  \sum_{q'}\int_{\mathbb{R}^n} \alpha_{t+1}^i(x',q') Q_{q',y^q}(u) \varphi(y^x - h( x',u)) \tau(x',q'\mid x,q,u)\,dx' {\bf 1}_K(x,q)
\end{equation}
Because $\alpha_{y,u}^i$ does not depend on $\sigma$, we can redefine the supremum over all $\Gamma_{t+1}$ to be over all $\alpha_{y,u}^i$.
\begin{equation*}
V_t^*(\sigma) =  \sup_{u\in\mathcal{U}} \sum_{y^q} \int_{\mathbb{R}^n} \sup_{\{\alpha_{y,u}^i\}} \,\langle \alpha_{y,u}^i, \sigma \rangle \,dy^x
\end{equation*}
For a specific $\sigma$, $u$, and $y$, if we define
\begin{equation}\label{alpha_yus}
\alpha_{y,u,\sigma}(x,q) = \arg\sup_i\, \langle \alpha_{y,u}^i, \sigma \rangle
\end{equation}
then $V_t^*$ can be further simplified as 
\begin{equation*}
V_t^*(\sigma) =  \sup_{u\in\mathcal{U}} \sum_{y^q} \int_{\mathbb{R}^n} \langle \alpha_{y,u,\sigma}, \sigma\rangle \,dy^x =  \sup_{u\in\mathcal{U}} \,\left\langle  \sum_{y^q} \int_{\mathbb{R}^n} \alpha_{y,u,\sigma}\,dy^x, \sigma \right\rangle
\end{equation*}
Therefore, the set of all $\{\alpha_t^i\}$ can be described by
\begin{equation}\label{alpha_set_cont}
\Gamma_t = \bigcup_{\sigma}\left\{\sum_{y^q} \int_{\mathbb{R}^n} \alpha_{y,u,\sigma}\,dy^x\right\}_{\forall u\in\mathcal{U}}
\end{equation}
and $V_t^*$ may be written as 
\begin{equation}
V_t^*(\sigma) = \sup_{\alpha_t^i\in\Gamma_t} \,\langle \alpha_t^i, \sigma \rangle
\end{equation}

\end{IEEEproof}

As in \cite{Porta06}, for discrete actions and observations, the set $\Gamma_t$ has finite cardinality, and so $V_t^*(\sigma)$ is a piecewise-linear function in $\sigma$.  If the state space was small and discrete, as were the observations and actions, we could construct a finite set of $\sigma$ vectors, and then generate a finite set of $\alpha$-vectors at each time step $k$ to solve the above problem exactly, much like the algorithm first proposed by \cite{sondik}.  However, with a hybrid state space, there are an infinite number of $\sigma$ functions defined on an infinite number of states, and so we cannot hope to solve this problem exactly.  We can instead sample sufficient statistics $\sigma$ from the set $\Sigma$ of all possible $\sigma$ functions, just as a collection of sampled belief points are used in \cite{Porta06} and many other PBVI solvers designed for large (but discrete) state spaces.  The set of sampled points is denoted $\tilde{\Sigma}$.   By sampling from the sufficient statistic space $\Sigma$, we can generate a finite number of $\alpha$-functions.  Further, because of the piecewise-linear convex nature of the value functions, we are guaranteed to obtain a lower bound on the true value function.  In fact, we can characterize the error between the value functions produced by the point-based method and the true value functions, based on how densely we sample $\Sigma$.

The operator $H$ in \eqref{DPop}  represents the complete backup operation \eqref{backupVal}.  The operator $\tilde{H}$ will be used to represent a point-based backup based on a set of sampled belief points $\tilde{\Sigma}$.  We denote the approximate value function at time $t$ characterized by $\tilde{\Gamma}_t$ as $V_t^{\tilde{\Sigma}}$, in comparison to the true value function $V_t^*$.  Further, let $\delta(\tilde{\Sigma})$ be the maximum $L^1$ distance between points in $\tilde{\Sigma}$ and points in $\Sigma$.
\begin{equation}\label{delta}
\delta(\tilde{\Sigma}) = \sup_{\tilde{\sigma}\in\tilde{\Sigma}}\inf_{\sigma\in\Sigma}\|\tilde{\sigma} - \sigma\|_1
\end{equation}
Now consider the maximum error introduced in performing one iteration of point-based backup, given the current value function estimate $V_t^{\tilde{\Sigma}}$.

\begin{lem}\label{lemma3}
The error introduced in one iteration of point-based value iteration, denoted $\epsilon^{(1)}$, is at most $\delta(\tilde{\Sigma})$:
\begin{equation*}
 \left\|\tilde{H}[V_t^{\tilde{\Sigma}}] - H[V_t^{\tilde{\Sigma}}]\right\|_{\infty} = \epsilon^{(1)} \leq \delta(\tilde{\Sigma})
\end{equation*}
\end{lem}
\begin{IEEEproof}
The proof is similar to one  in \cite{pineau06}  for discrete state POMDPs.  First, let $\sigma^{(1)}$ be the point in $\Sigma$ where the error between the true value function and the point-based backup is greatest.  Let $\sigma^{(2)} \in \tilde{\Sigma}$ be the closest point in the $L^1$ sense to $\sigma^{(1)}$.  Let $\alpha^{(2)} \in \tilde{\Gamma}_{t-1}$ be maximal at $\sigma^{(2)}$, and $\alpha^{(1)}\in \Gamma_{t-1}$ (and not in $\tilde{\Gamma}_{t-1}$) is the function that \emph{would} be maximal at $\sigma^{(1)}$ had it been calculated.

Then
\begin{align}
\epsilon^{(1)} &\leq |\langle\alpha^{(1)}, \sigma^{(1)}\rangle - \langle \alpha^{(2)}, \sigma^{(1)}\rangle| \notag\\
&\leq |\langle\alpha^{(1)}, \sigma^{(1)}\rangle - \langle \alpha^{(2)}, \sigma^{(1)}\rangle + \langle\alpha^{(1)}, \sigma^{(2)}\rangle - \langle \alpha^{(1)}, \sigma^{(2)}\rangle| \notag \\
& \leq |\langle\alpha^{(1)}, \sigma^{(1)}\rangle - \langle \alpha^{(2)}, \sigma^{(1)}\rangle + \langle\alpha^{(2)}, \sigma^{(2)}\rangle - \langle \alpha^{(1)}, \sigma^{(2)}\rangle| \label{VoptProp}\\
&\leq |\langle \alpha^{(1)} - \alpha^{(2)}, \sigma^{(1)} - \sigma^{(2)}\rangle| \notag 
\end{align}
\begin{align}
& \leq \|\alpha^{(1)} - \alpha^{(2)}\|_{\infty} \| \sigma^{(1)} - \sigma^{(2)}\|_1  \label{Holder}\\
&\leq  \|\alpha^{(1)} - \alpha^{(2)}\|_{\infty}  \delta(\tilde{\Sigma})  \label{last}
\end{align}
Line \eqref{VoptProp} follows because $\alpha^{(2)}$ is optimal for $\sigma^{(2)}$, implying $\langle\alpha^{(1)}, \sigma^{(2)}\rangle \leq \langle \alpha^{(2)}, \sigma^{(2)}\rangle$.  Line \eqref{Holder} follows from H\"{o}lder's Inequality.  Line \eqref{last} can be further simplified by noting that the $\alpha$-functions are bounded between 0 and 1 for all $x\in\mathcal{X}$ and $q\in\mathcal{Q}$.  Because the value function at a specific point $\sigma$ represents the probability of staying within set $K$ for some length of time, given the normalized density $\sigma$, this value must be between 0 and 1.  The value function is further defined as $\sup\langle \alpha, \sigma\rangle$, meaning that the inner product of $\alpha$ and $\sigma$ must be between 0 and 1, and therefore $\alpha$ must be between 0 and 1 for all $x,q$ (since by \eqref{alpha_yu} it clearly must be nonnegative).

Therefore, we can say $\| \alpha^{(1)} - \alpha^{(2)}\|_{\infty}\leq 1$, and we get that $\epsilon^{(1)} \leq \delta(\tilde{\Sigma})$.

\end{IEEEproof}

We now use Lemma \ref{lemma3} to derive a bound between the true value function and the point-based approximation at any time $t$.

\begin{thm}\label{thm1}
For a set of sufficient statistics $\Sigma$, sampled set $\tilde{\Sigma}$, and horizon $t$, the error from using point-based value iteration versus full value iteration, given by $\epsilon(t) = \|V_t^{\tilde{\Sigma}} - V_t^*\|_{\infty}$ is bounded above by
\begin{equation*}
 \|V_t^{\tilde{\Sigma}} - V_t^*\|_{\infty}= \epsilon(t) \leq t \delta(\tilde{\Sigma})
\end{equation*}
\end{thm}
\begin{IEEEproof}
\begin{align}
\epsilon(t) &= \|V_{T-t}^{\tilde{\Sigma}} - V_{T-t}^*\|_{\infty}\notag\\
&= \|\tilde{H}[V_{T-t-1}^{\tilde{\Sigma}}] - H[V_{T-t-1}^*]\|_{\infty} \notag\\
&=  \|\tilde{H}[V_{T-t-1}^{\tilde{\Sigma}}] - H[V_{T-t-1}^*] + H[V_{T-t-1}^{\tilde{\Sigma}}] -  H[V_{T-t-1}^{\tilde{\Sigma}}] \|_{\infty}\notag\\
&\leq \| \tilde{H}[V_{T-t-1}^{\tilde{\Sigma}}]-  H[V_{T-t-1}^{\tilde{\Sigma}}]\|_{\infty} + \| H[V_{T-t-1}^{\tilde{\Sigma}}] -  H[V_{T-t-1}^*] \|_{\infty} \notag\\
&\leq \epsilon^{(1)} + \|V_{T-t-1} ^{\tilde{\Sigma}} - V_{T-t-1}^*\|_{\infty} \label{bydef}\\
&\leq \epsilon^{(1)} + \epsilon(t-1) \notag \\
\epsilon(t) &\leq t\delta(\tilde{\Sigma}) \label{bylem}
\end{align}
Line \eqref{bydef} follows from the definition of $\epsilon^{(1)}$, and line \eqref{bylem} follows from Lemma \ref{lemma3}.  
\end{IEEEproof}

Thus the error between the point-based approximation and the actual value function is directly proportional to how densely $\tilde{\Sigma}$ is sampled, and converges to zero as $\tilde{\Sigma}$ approaches $\Sigma$.

\subsection{Implementation}

For a state space $\mathcal{S}$ that is discrete, ``closedness'' of the belief state $b(s)$ and of the $\alpha$-vectors $\alpha_t(s)$ is maintained after updates by the operator $M_{y,u}$ and by \eqref{alpha1}-\eqref{alpha2}, respectively.  That is, although the belief function can take on an infinite number of values for each state $s$ (the interval $[0,1]$), because there are a finite number of states in $\mathcal{S}$, the function $b(s)$ can be represented by a vector $[b(s_0) \, b(s_1) \dotsc b(s_n)]$ with each entry $b(s_i)\in [0,1]$ corresponding to the probability of being in state $s_i$ according to the specific density $b$.  Similarly for the set of $\alpha$-vectors, $\Gamma_t$, which remain the same size after updates according to \eqref{alpha1} and \eqref{alpha2}.  

For $\mathcal{S}$ continuous, this ``closedness'' property of the structure of both the beliefs and $\alpha$-functions under updating is no longer guaranteed, and can make the computation intractable.  As a remedy, \cite{Porta06} represents both the beliefs and $\alpha$-functions as sums of weighted Gaussians (which can represent a function to any desired accuracy with enough components), and shows that for an additive cost POMDP, the belief function remains a Gaussian sum under the belief update operator $M_{y,u}$, as do the $\alpha$-functions when generated recursively from the previous set of $\alpha$-functions.  The Gaussian sum representation also guarantees the inner product operation $\langle \alpha, b\rangle$ to be computable. 

We now show that we can approximate the sufficient statistic $\sigma$ by a vector whose entries are finite sums of Gaussians (each entry of the vector corresponds to a different discrete mode $q$), and that this representation is closed under the update operator $\Phi$.  We also show that the $\alpha$-functions as defined by \eqref{alpha_set_cont} for the multiplicative reachability cost function can also be approximated by vectors of finite sums of Gaussians, and are closed under the operations defined in \eqref{alpha_yu} and \eqref{alpha_yus}.  All of the following derivations assume a discrete observation space of finite cardinality $N_{y^q}\times N_{y^x}$.  We make the following additional assumptions:

\emph{Assumption 1:} We can represent the indicator function as a finite sum of Gaussians \eqref{IndApprox}, with $w_i(q)\in\mathbb{R}$ a mode-dependent coefficient, such that for $q\in K_q$, $w_i(q)=1$, and for $q\notin K_q$ $w_i(q)=0$ for all $i$, where $K=K_x\times K_q$.  Gaussian distribution $i$ has mean $\mu_i$ and covariance $\Sigma_i$.
\begin{equation}\label{IndApprox}
{\bf 1}_K(x,q) \approx \sum_{i=1}^I w_i(q)\phi(x; \mu_i, \Sigma_i)
\end{equation}

\emph{Assumption 2:} We can approximate the stochastic kernel $\tau(s'\mid s, u) = T_x(x'\mid x, q',u)T_q(q'\mid x, q, u)$ by a Gaussian sum.  We first express the distribution of the discrete variable $q'$ in terms of Gaussian distributions evaluated at the continuous variable $x$:
\begin{equation}\label{Tqsum}
T_q(q'\mid x, q, u) \approx \sum_{j=1}^J w_j(q',q,u) \phi(x; \mu_j(q',q,u), \Sigma_j(q',q,u))
\end{equation}
For finite $J$ \eqref{Tqsum} will never exactly sum to 1 (see \cite{Brunskill2010}) and so will always be an approximation.  We assume that the continuous dynamics are linear in $x$ with Gaussian noise, so that
\begin{equation}\label{Txsum}
T_x(x'\mid x, q', u)= \phi(x'; \mu_{q'}^u(x), \mathcal{W}_{q'}^{u})
\end{equation} 
where $\mu_{q'}^u(x)$ is of the form $Ax+f(q',u)$ with $A\in\mathbb{R}^{n\times n}$ invertible and $f$ a possibly non-linear function of $q'$ and $u$.  This allows us to rewrite $T_x$ in terms of $x$ rather than $x'$, so that $T_x(x'\mid x,q',u) = \delta\phi(x; \hat{\mu}_{q'}^u(x'), \hat{\mathcal{W}}_{q'}^u)$ as well.  In fact, 
\begin{equation}\label{gauss_transform}
\phi(x'; Ax+f(q',u), \mathcal{W}_{q'}^u) = |A^{-1}|\phi\left(x; A^{-1}(x'-f(q',u)), A^{-1}\mathcal{W}_{q'}^u(A^{-1})^T\right)
\end{equation}
with $\hat{\mu}_{q'}^u(x') = A^{-1}(x'-f(q',u))$ and $\hat{\mathcal{W}}_{q'}^u = A^{-1}\mathcal{W}_{q'}^u(A^{-1})^T$.

\emph{Assumption 3:} The discrete observation model for the continuous variable, $\varphi(y^x\mid x,u)$ can be approximated by
\begin{equation}\label{ysum}
\varphi(y^x\mid x,u) \approx \sum_{h=1}^Hw_h(y,u) \phi(x; \mu_h(y,u), \Sigma_h(y,u))
\end{equation}

To make notation (slightly) cleaner, we now shift any parameter's dependence on either $y$ or $u$ to its superscript, and any dependence on $q$ or $q'$ to its subscript, so for instance $w_h(y,u)$ becomes $w_h^{y,u}$ and $\mu_j(q',q,u)$ becomes $\mu_{j,q',q}^{u}$. 

\subsection{Approximating the Sufficient Statistic} \label{sigma_derivation}

\begin{lem}\label{lem:sigma}
The sufficient statistic $\sigma_t(x,q)$ can be approximated by a linear combination of Gaussians for all $t= 0,1,\dotsc$, where the parameters of each Gaussian component are dependent on the discrete variable $q$.
\begin{equation}\label{sigmaform}
\sigma_t(x,q) \approx \sum_{l=1}^Lw_{l,q}\phi(x; \mu_{l,q},\Sigma_{l,q})
\end{equation}  
\end{lem}
\begin{IEEEproof}
The proof follows by induction.  For $t=0$, $\sigma_0(x,q) = \rho(x,q)$.  Because any distribution can be approximated to arbitrary accuracy by a weighted sum of Gaussians, we set 
$\rho(x,q) = \sum_{l=1}^L w_{l,q}\phi(x; \mu_{l,q}, \Sigma_{l,q})$  and so $\sigma_0(x,q)$ is of the form \eqref{sigmaform}.

For $t=n-1$, assume that $\sigma_{n-1}(x,q) =   \sum_{l=1}^Lw_{l,q}\phi(x; \mu_{l,q},\Sigma_{l,q})$.  Then under the operator $\Phi$ given by \eqref{sigma_op}, it follows that

\begin{align*}
\sigma_n(x',q') &= N_{y^q}N_{y^x}Q_{q',y^q}(u) \psi(y^x\mid x', u) \sum_{q=1}^{N_q}\int_{\mathbb{R}^n} {\bf 1}_K(x,q) T_x(x'\mid x,q',u)T_q(q'\mid q, x, u) \sigma_{n-1}(x,q)\,dx \\
&\approx N_{y^q}N_{y^x}Q_{q',y^q}(u)\left[\sum_{h=1}^Hw_h^{y,u}\phi(x';\mu_h^{y,u}, \Sigma_h^{y,u})\right] \sum_{q=1}^{N_q}\int_{\mathbb{R}^n} \left[\sum_{i=1}^Iw_{i,q}\phi(x; \mu_{i}, \Sigma_{i})\right] \\
&\hspace{15 mm}\times |A^{-1}|\phi\left(x; A^{-1}(x'-f(q',u)), A^{-1}\mathcal{W}_{q'}^u(A^{-1})^T\right) \left[\sum_{j=1}^J w_{j,q,q'}^{u}\phi(x; \mu_{j,q,q'}^{u}, \Sigma_{j,q,q'}^{u})\right]\\
&\hspace{15 mm} \times \left[\sum_{l=1}^Lw_{l,q}\phi(x; \mu_{l,q}, \Sigma_{l,q})\right]\,dx\\
&\approx\sum_{h=1}^H\sum_{i=1}^I\sum_{j=1}^J\sum_{l=1}^L\sum_{q=1}^{N_q}N_{y^q}N_{y^x}Q_{q',y^q} (u)|A^{-1}|w_h^{y,u}w_{i,q}w_{j,q,q'}^uw_{l,q}\phi(x'; \mu_h^{y,u}, \Sigma_h^{y,u}) \\
&\hspace{15 mm} \times \int_{\mathbb{R}^n} \phi(x; \mu_{i}, \Sigma_{i})\phi\left(x; \hat{\mu}_{q'}^u(x'), \hat{\mathcal{W}}_{q'}^u\right)\phi(x; \mu_{j,q,q'}^{u}, \Sigma_{j,q,q'}^{u}) \phi(x; \mu_{l,q}, \Sigma_{l,q}) \,dx \\
\end{align*}
Next, the below identity regarding multiplied Gaussians is used to combine the above Gaussians inside the integral.
\begin{equation} \label{Gaussprod} 
\begin{split}
\phi(x;\mu_1, \Sigma_1)\phi(x; \mu_2, \Sigma_2) &= \phi(\mu_1; \mu_2, \Sigma_1+\Sigma_2)\phi(x; \tilde{\mu},\tilde{\Sigma})\\
\tilde{\mu} &= \tilde{\Sigma}(\Sigma_1^{-1}\mu_1 + \Sigma_2^{-1}\mu_2)   \\
\tilde{\Sigma} &= (\Sigma_1^{-1} + \Sigma_2^{-1})^{-1} 
\end{split}
\end{equation}
Then 
\begin{align*}
\sigma_n(x',q') &\approx \sum_{h,i,j,l,q} N_{y^q}N_{y^x}Q_{q',y^q} (u)|A^{-1}|w_h^{y,u}w_{i,q}w_{j,q,q'}^uw_{l,q}\phi(x'; \mu_h^{y,u}, \Sigma_h^{y,u})\phi(\hat{\mu}_{q'}^u(x'); \mu_{i}, \hat{\mathcal{W}}_{q'}^u+\Sigma_{i}) \\
&\hspace{30 mm} \times \phi(\mu_{j,q,q'}^u; \mu_{l,q}, \Sigma_{j,q,q'}^u + \Sigma_{l,q}) \int_{\mathbb{R}^n} \phi(x; \tilde{\mu}_1(x'), \tilde{\Sigma}_1)\phi(x; \tilde{\mu}_2, \tilde{\Sigma}_2)\,dx \\ 
\end{align*}
with
\begin{align*}
\tilde{\mu}_1(x') &= \tilde{\Sigma}_1\left(\Sigma_{i}^{-1}\mu_{i} +\left(\hat{\mathcal{W}}_{q'}^u\right)^{-1}\hat{\mu}_{q'}^u(x')\right), \hspace{5 mm}
\tilde{\Sigma}_1 = \left(\Sigma_{i}^{-1} +\left(\hat{\mathcal{W}}_{q'}^u\right)^{-1}\right)^{-1} \\
\tilde{\mu}_2 &= \tilde{\Sigma}_2 \left(\left(\Sigma_{j,q,q'}^u\right)^{-1}\mu_{j,q,q'}^u + \Sigma_{l,q}^{-1}\mu_{l,q}\right), \hspace{5 mm}
\tilde{\Sigma}_2 = \left(\left(\Sigma_{j,q,q'}^u\right)^{-1}+\Sigma_{l,q}^{-1}\right)^{-1}
\end{align*}
using \eqref{Gaussprod}.  Multiplying the final two Gaussians inside the integral leaves only one Gaussian that is a function of $x$, which integrates to $1$, leaving
\begin{align*}
\sigma_n(x,q) &\approx \sum_{h,i,j,l,q}N_{y^q}N_{y^x}Q_{q',y^q} (u)|A^{-1}|w_h^{y,u}w_{i,q}w_{j,q,q'}^uw_{l,q}\phi(\mu_{j,q,q'}^u; \mu_{l,q}, \Sigma_{j,q,q'}^u + \Sigma_{l,q})\phi(x'; \mu_h^{y,u}, \Sigma_h^{y,u}) \\
&\hspace{30 mm} \times \phi(\hat{\mu}_{q'}^u(x'); \mu_{i}, \hat{\mathcal{W}}_{q'}^u+\Sigma_{i})\phi(\tilde{\mu}_1(x'); \tilde{\mu}_2, \tilde{\Sigma}_1 + \tilde{\Sigma}_2)
\end{align*}
Now all that is left to complete the proof is to manipulate the last two Gaussians, $\phi(\hat{\mu}_{q'}^u(x'); \mu_{i,q}, \newline\hat{\mathcal{W}}_{q'}^u +\Sigma_{i})$ and $\phi(\tilde{\mu}_1(x'); \tilde{\mu}_2, \tilde{\Sigma}_1 + \tilde{\Sigma}_2)$ so that they are functions of $x'$, i.e. $\phi(x'; \hat{\mu}, \hat{\Sigma})$, and then apply \eqref{Gaussprod} twice.  This can be done in both cases using straightforward but tedious linear algebra. First noting that
\begin{align*}
\phi(\hat{\mu}_{q'}^u(x'); \mu_{i}, \hat{\mathcal{W}}_{q'}^u+\Sigma_{i}) & =|A|\phi(x';A\mu_{i}+f(q',u), \mathcal{W}_{q'}^u+A\Sigma_{i}A^T)\\
\phi(\tilde{\mu}_1(x'); \tilde{\mu}_2, \tilde{\Sigma}_1 + \tilde{\Sigma}_2) &=|A\hat{\mathcal{W}}_{q'}^u\tilde{\Sigma}_1^{-1}|\phi(x';\overline{\mu}_{q,q'}^u, \overline{\Sigma}_{q,q'}^u)\\
\overline{\mu}_{q,q'}^u &=A\left[\hat{\mathcal{W}}_{q'}^u\tilde{\Sigma}_1^{-1}(\tilde{\mu}_2-\tilde{\Sigma}_1\Sigma_{i}^{-1}\mu_{i})+f(q',u)\right] \\
\overline{\Sigma}_{q,q'}^u &= A\left[\hat{\mathcal{W}}_{q'}^u\tilde{\Sigma}_1^{-1}(\tilde{\Sigma}_1+\tilde{\Sigma}_2)\hat{\mathcal{W}}_{q'}^u\tilde{\Sigma}_1^{-1}\right]A^T
\end{align*}
we can ultimately write
\begin{equation}\label{sig1}
\begin{split}
\sigma_n(x',q')&\approx \sum_{h,i,j,l,q}w_{h,i,j,l,q,q'}\phi(x'; \mu_{h,i,j,l,q,q'}, \Sigma_{h,i,j,l,q,q'}) \\
&\approx \sum_{k=1}^{HIJLN_q}w_{k,q'}\phi(x'; \mu_{k,q'}, \Sigma_{k,q'}) 
\end{split}
\end{equation}
where
\begin{align}
\begin{split}\label{sig2}
w_{h,i,j,l,q,q'} &= |A\hat{\mathcal{W}}_{q'}^u\tilde{\Sigma}_1^{-1}| N_{y^q}N_{y^x}Q_{q',y^q} (u)w_h^{y,u}w_{i,q}w_{j,q,q'}^uw_{l,q}\phi(\mu_{j,q,q'}^u; \mu_{l,q}, \Sigma_{j,q,q'}^u + \Sigma_{l,q}) \\
&\hspace{20 mm} \times \phi(\mu_h^{y,u}; A\mu_{i}+f(q',u), \Sigma_h^{y,u}+\mathcal{W}_{q'}^u+A\Sigma_{i}A^T)\phi(\overline{\mu}_{q,q'}^u; c, C + \overline{\Sigma}_{q,q'}^u)
\end{split}\\
\mu_{h,i,j,l,q,q'} &= \left(C^{-1}+\left(\overline{\Sigma}_{q,q'}^u\right)^{-1}\right)^{-1}\left(C^{-1}c + \left(\overline{\Sigma}_{q,q'}^u\right)^{-1}\overline{\mu}_{q,q'}^u\right) \\
\Sigma_{h,i,j,l,q,q'} &=  \left(C^{-1}+\left(\overline{\Sigma}_{q,q'}^u\right)^{-1}\right)^{-1}
\end{align}
\begin{align}
C &= \left(\left(\Sigma_h^{y,u}\right)^{-1}+\left(\mathcal{W}_{q'}^u+A\Sigma_{i}A^T\right)^{-1}\right)^{-1}\\
c &= C\left(\left(\Sigma_h^{y,u}\right)^{-1}\mu_h^{y,u}  + \left(\mathcal{W}_{q'}^u+A\Sigma_{i}A^T\right)^{-1}\left(A\mu_{i}+f(q',u)\right)\right) \label{siglast}
\end{align}
\end{IEEEproof}
The  sufficient statistic $\sigma$ is therefore closed under the update operator $\Phi$.  The expression in \eqref{sig1} - \eqref{siglast} simplifies somewhat depending on the problem, as seen in Section \ref{examples}.  More problematic is the explosion in the number of Gaussians: for $L$ Gaussians representing $\sigma_{n-1}$, $HIJLN_{y^q}$ Gaussians are required to represent $\sigma_n$.  However, there are techniques to combine similar components (the individual weighted Gaussians) of the mixture in order to bound the total components, which will be discussed in Section \ref{examples}.

\subsection{Approximating the $\alpha$-Functions} \label{alpha_derivation}

We use the same approach as in Lemma \ref{lem:sigma} to approximate the $\alpha$-functions by Gaussian mixtures, through induction and application of the operation defined in \eqref{alpha_yu} that generates $\alpha_{y,u}^i$ from $\alpha_{n+1}^i$.  Showing that \eqref{alpha_yu} preserves the Gaussian mixture structure of the $\alpha$-functions is sufficient to show that the full $backup$ operation is closed under Gaussian sums when the observations $y^x$ are discrete, since the only additional operation is to sum over all $y^x$, as in \eqref{alpha_set_cont}.

\begin{lem}\label{lem:alpha}
The $\alpha$-functions $\alpha_t^i(x,q)$ can be approximated by a linear combination of Gaussians for all $t=0,1,\dotsc$ where the parameters for each Gaussian are dependent on the discrete variable $q$.
\begin{equation}\label{alpha_approx}
\alpha_t^i(x,q) \approx \sum_{d=1}^D w_{d,q}\phi(x;\mu_{d,q}, \Sigma_{d,q})
\end{equation}
\end{lem}
\begin{IEEEproof}
We omit most details of the proof, since they are almost identical to those in the proof of Lemma \ref{lem:sigma}. 

 For $t=T$, from Lemma \ref{lemma2} and the definition of $\alpha_T$ as the indicator function ${\bf 1}_K(s)$, setting
\begin{equation*}
\alpha_T(x,q) \approx \sum_{i=1}^I w_{i,q}\phi(x; \mu_i, \Sigma_i)
\end{equation*}
using the Gaussian sum approximation to the indicator function \eqref{IndApprox} gives $\alpha_T$ in the desired form.

Assuming $\alpha_{n+1}^j(x',q') = \sum_{d=1}^{D}w_{d,q'}\phi(x'; \mu_{d,q'}, \Sigma_{d,q'})$, using \eqref{alpha_yu} it follows that

\begin{align*}
\alpha_{y,u}^j(x,q) &\approx \sum_{q'=1}^{N_q}Q_{q',y^q}(u)\int_{\mathbb{R}^n} \left[\sum_{d=1}^Dw_{d,q'}\phi(x'; \mu_{d,q'}, \Sigma_{d,q'})\right]\left[\sum_{h=1}^Hw_h^{y,u}\phi(x';\mu_h^{y,u},\Sigma_h^{y,u})\right] \\
&\hspace{10 mm}\times \phi(x';\mu_{q'}^u(x), \mathcal{W}_{q'}^u)\,dx'\left[\sum_{j=1}^Jw_{j,q',q}^u\phi(x;\mu_{j,q',q}^u,\Sigma_{j,q',q}^u)\right]\left[\sum_{i=1}^Iw_{i,q}\phi(x;\mu_i,\Sigma_i)\right] \\
&\approx \sum_{q',d,h,j,i} w_{q',d,h,j,i,q}\phi(x; \mu_{q',d,h,j,i,q}, \Sigma_{q',d,h,j,i,q})
\end{align*}
where
\begin{align}
\begin{split}
w_{q',d,h,j,i,q} &= |A^{-1}|Q_{q',y^q}(u)w_{d,q'}w_h^{y,u}w_{j,q',q}^uw_{i,q}\phi(\mu_i; \mu_{j,q',q}^u,\Sigma_i+\Sigma_{j,q',q}^u)\phi(\mu_{d,q'}; \mu_h^{y,u}, \Sigma_{d,q'}+\Sigma_h^{y,u})  \\
&\hspace{20 mm}\times \phi\left(\tilde{\mu}_1; A^{-1}\left(\tilde{\mu}_2-f(q',u)\right), \tilde{\Sigma}_1+A^{-1}\left(\tilde{\Sigma}_2+\mathcal{W}_{q'}^u\right)\left(A^{-1}\right)^T\right) 
\end{split}\\
\mu_{q',d,h,j,i,q} &= \left(A^T\left(\tilde{\Sigma}_2+\mathcal{W}_{q'}^u\right)^{-1}A+\tilde{\Sigma}_1^{-1}\right)^{-1}  \left(A^T\left(\tilde{\Sigma}_2+\mathcal{W}_{q'}^u\right)^{-1}\left(\tilde{\mu}_2-f(q',u)\right)+\tilde{\Sigma}_1^{-1}\tilde{\mu}_1\right) \\
\Sigma_{q',d,h,j,i,q} &= \left(A^T\left(\tilde{\Sigma}_2+\mathcal{W}_{q'}^u\right)^{-1}A+\tilde{\Sigma}_1^{-1}\right)^{-1} \\
\tilde{\mu}_1 &= \tilde{\Sigma}_1\left(\left(\Sigma_{j,q',q}^u\right)^{-1}\mu_{j,q',q}^u + \Sigma_i^{-1}\mu_i\right), \hspace{5 mm}
\tilde{\Sigma}_1 = \left(\left(\Sigma_{j,q',q}^u\right)^{-1} + \Sigma_i^{-1}\right)^{-1} \notag\\
\tilde{\mu}_2 &= \tilde{\Sigma}_2\left(\Sigma_{d,q'}^{-1}\mu_{d,q'}+\left(\Sigma_h^{y,u}\right)^{-1}\mu_h^{y,u}\right), \hspace{5 mm}
\tilde{\Sigma}_2 =\left(\Sigma_{d,q'}^{-1} + \left(\Sigma_h^{y,u}\right)^{-1}\right)^{-1}\notag
\end{align}
Because each $\alpha_{y,u}^j(x,q)$ is approximated by a sum of Gaussians for all $j$, any $\alpha_n^j\in \tilde{\Gamma}_n$ is also a sum of Gaussians (since the only additional operation to generate the $\alpha_n^j$ from $\alpha_{y,u}^j$ is to sum over all $y^x$ and $y^q$, as in \eqref{alpha_set_cont}).
\end{IEEEproof}

\section{Example}\label{examples}

The temperature regulation problem is a benchmark example for hybrid systems, and a stochastic version with perfect state information is presented in \cite{Abatecomp}.  We consider the case of one heater, which can either be turned on to heat one room, or turned off.  The temperature of the room at time $t$ is given by the continuous variable $x(t)$, and the discrete state $q(t)=1$ indicates the heater is on at time $t$, and $q(t)=0$ denotes the heater is off.  The stochastic difference equation governing the temperature is given by
\begin{equation*}
x(t+1) = (1-b)x(t)+cq(t+1)+bx_a+v(t)
\end{equation*}
with constants $b=0.0167$, $c=0.8$, and $x_a=6$, and $v(t)$  i.i.d. Gaussian random variables with mean zero and variance $v^2$.  The control input is given by $u(t)\in\mathcal{U}$ with $\mathcal{U}=\{0,1\}$, but the chosen control is not always implemented with probability $1$.  Instead, $q(t)$ is updated probabilistically, dependent on $u(t-1)$ and $q(t-1)$, with transition function $T_q(q(t+1)\mid q(t), u(t))$.  So while function $\overline{\mu}_t(\sigma_t)$ deterministically returns a single control input, control input $u_t = \overline{\mu}_t(\sigma_t)$ may not always be implemented.

To model this as a partially observable problem, assume the actual temperature is unknown, and only a noisy measurement is available to the controller.  The controller does, however, know if the heater is on or off at time $t$ (i.e. $q(t)$ is perfectly observed).  The observation $y(t) = y^x(t)$  is given by $y^x(t) = x(t) + w(t)$, with $w(t)$ i.i.d. Gaussian random variables with mean zero and variance $w^2$ (so that $\varphi(\hat{w}) = \phi(\hat{w}; 0, w^2)$.  Because the discrete mode $q$ is perfectly observed, we do not record $y^q(t)$, and it is not included in any equations.

It is desirable to keep the temperature of the room between $17.5$ and $22$ degrees celsius at all times, hence the safe region $K=[17.5,22]$ does not depend on the discrete state $q(k)$ (so ${\bf 1}_K(s) = {\bf 1}_K(x)$).  To find the maximum probability that the room stays within the desired temperature range given that the controller only has access to the mode $q(t)$ and observations $y^x(t)$, we first find expressions for both $\sigma_t(x,q)$ and $\alpha_t(x,q)$ as Gaussian sums.

We first discretize the observations $y^x$.  Using a grid with spacing $\Delta$, $y^x$ is redefined over $\mathcal{Y} = \{17.5-tol_1, 17.5-tol_1+\Delta,\dotsc,22+tol_2-\Delta, 22+tol_2\}$ where $tol_1$ and $tol_2$ are defined so that the probability of observing $y^x$ outside of $[17.5-tol_1,22+tol_2]$ is approximately zero.  The probability that $y^x = \overline{y}\in\mathcal{Y}$ can be written as
\begin{equation*}\label{y_approx}
\varphi(\overline{y}-x)= \mathbb{P}\left[y^x\in\left[\left. \overline{y}-\frac{\Delta}{2}, \overline{y}+\frac{\Delta}{2}\right]\right| x\right] = \int_{\overline{y}-\frac{\Delta}{2}}^{ \overline{y}+\frac{\Delta}{2}} \phi(y; x, w^2)\, dy
\end{equation*}
which in turn can be approximated by a summation:
\begin{equation*}
\int_{\overline{y}-\frac{\Delta}{2}}^{ \overline{y}+\frac{\Delta}{2}} \phi(y; x, w^2)\, dy = \int_{\overline{y}-\frac{\Delta}{2}}^{ \overline{y}+\frac{\Delta}{2}} \phi(x; y, w^2)\, dy
\approx \sum_{h=\overline{y}-\frac{\Delta}{2}}^{ \overline{y}+\frac{\Delta}{2}-\delta_H} \delta_H \phi(x; h, w^2)
\end{equation*}
For $\delta_H$ the grid spacing in the interval $\left[\overline{y}-\frac{\Delta}{2}, \overline{y}+\frac{\Delta}{2}\right]$, we can now write the discretized observation function as a sum of Gaussians:
\begin{equation}\label{obs_approx}
\varphi(y^x-x) \approx \sum_{h=1}^H w_h \phi(x; \mu_h(y^x), w^2)
\end{equation}
with $w_h=\delta_H$ for all $h$, and $\mu_h(y^x) = y^x-\frac{\Delta}{2} + (h-1)\delta_H$.

The one dimensional indicator function also needs to be approximated by a sum of Gaussians.  Unfortunately, because the indicator function is discontinuous, approximation by a finite sum of Gaussians induces a pseudo-Gibbs phenomenon, with oscillations occuring near the discontinuities (endpoints of $K$).  Using more components leads to a smoother approximation in the interior of $K$, but the oscillations at the endpoints remain.  Unfortunately, no clear ways to avoid this phenomenon currently exist.  It should be noted that the inability to exactly represent the indicator function using a Gaussian mixture leads to  $\alpha$-functions that are also only approximations to the true $\alpha$-functions, and thus the guaranteed lower bound on the value function breaks down. One practical workaround is to choose Gaussian components that slightly underapproximate the indicator function (except at the endpoints), to help preserve the underapproximation to the true value function.  For low-dimensional problems we have not experienced any problems approximating the indicator function and obtaining a reasonable lower bound for the value function, but at higher dimensions it is possible that the number of Gaussian components required for reasonable approximations becomes prohibitive.

A recursive expression for $\sigma_k(x,q)$ can be found using the derivation given in section \ref{sigma_derivation}.  For an initial distribution $\rho(x)$ on $x_0$ that is Gaussian with mean $\mu_0$ and variance $s^2$, and for $q_0=0$, then 
\begin{equation}\label{sig0}
\sigma_0(x,q) = {\bf 1}_{\{0\}}(q)\rho(x) = \begin{bmatrix} \sigma_0(x,0) \\ \sigma_0(x,1) \end{bmatrix} = \begin{bmatrix} \phi(x; \mu_0,s^2) \\ 0 \end{bmatrix}
\end{equation}

Given we already have an approximation to $\sigma_t(x,q)$ as in \eqref{sigmaform}, $\sigma_{t+1}$ corresponding to observation $y$ and control input $u$ can be written as
\begin{equation}\label{tempSigma}
\sigma_{t+1}(x,q) = \sum_{q_t=1}^{N_q}\sum_{i=1}^I\sum_{h=1}^H\sum_{l=1}^Lw_{i,h,l,q_t}(q,u,y)\phi\left(x; \mu_{i,h,l,q_t}(q,u,y), \Sigma_{i,h,l}\right)
\end{equation}
with
\begin{align*}
w_{i,h,l,q_t}(q,u,y) &= N_{y^x}T_q(q\mid q_k,u)w_iw_{l,q_t}w_h\phi(\mu_i; \mu_{l,q_t}, \Sigma_l+\Sigma_i)\\
&\hspace{45 mm}\times \phi(\mu_h(y); cq+bx_a+(1-b)\hat{\mu}, w^2 + \hat{\Sigma}(1-b)^2+v^2) \\
\mu_{i,h,l,q_t}(q,y) &= \frac{\mu_h(y)\left(\hat{\Sigma}(1-b)^2+v^2\right) + \left(cq+bx_a+(1-b)\hat{\mu}\right)w^2}{w^2 + v^2 + \hat{\Sigma}(1-b)^2}\\
\end{align*}
\begin{align*}
\Sigma_{i,h,l} &= \frac{w^2\left(\hat{\Sigma}(1-b)^2+v^2\right)}{w^2+v^2+\hat{\Sigma}(1-b)^2} \\
\hat{\mu} &=\frac{ \mu_i\Sigma_l+ \mu_{l,q_t}\Sigma_i}{\Sigma_l + \Sigma_i} , \hspace{15 mm} \hat{\Sigma} = \frac{\Sigma_l\Sigma_i}{\Sigma_l+\Sigma_i}
\end{align*}

Likewise, given $\alpha_{t+1}^j(x,q)$ in the form \eqref{alpha_approx}, $\alpha_{y,u}^j(x,q)$ corresponding to observation $y$ and control $u$ can be written
\begin{equation*}
\alpha_{y,u}^j(x,q) = \sum_{q_{t+1}=1}^{N_q}\sum_{i=1}^I\sum_{h=1}^H\sum_{l=1}^Lw_{i,h,l,q_{t+1}}(q,u,y)\phi\left(x; \mu_{i,h,l,q_{t+1}}(y),\Sigma_{i,h,l,q_{t+1}}\right)
\end{equation*}
with
\begin{align*}
w_{i,h,l,q_{t+1}}(q,u,y) & = \frac{1}{1-b}T_q(q_{t+1}\mid q, u)w_iw_{l,q_{t+1}}w_h\phi(\mu_h(y); \mu_{l},\Sigma_l+w^2) \\
&\hspace{50 mm} \times \phi\left(\mu_i; \frac{\hat{\mu}-cq_{t+1}-bx_a}{1-b}, \Sigma_i+\frac{\hat{\Sigma}+v^2}{(1-b)^2}\right) \\
\mu_{i,h,l,q_{t+1}}(y) &= \frac{\mu_i(\hat{\Sigma}+v^2)+(1-b)\Sigma_i(\hat{\mu}-cq_{t+1}-bx_a)}{\hat{\Sigma}+v^2+(1-b)^2\Sigma_i} \\
\Sigma_{i,h,l,q_{t+1}} &= \frac{\Sigma_i(v^2+\hat{\Sigma})}{\hat{\Sigma}+v^2+(1-b)^2\Sigma_i} \\
\hat{\mu} &= \frac{\mu_h(y)\Sigma_l+\mu_lw^2}{w^2+\Sigma_l}, \hspace{15 mm} \hat{\Sigma} = \frac{w^2\Sigma_l}{w^2 + \Sigma_l} 
\end{align*}

We implemented an algorithm in the style of POMDP solver Perseus \cite{Porta06} to generate an approximation to the value functions, by generating a fixed set of belief points, which we backed up at each iteration.   Unlike Perseus, we backed up every belief point, as necessary for a finite horizon calculation.  To generate the set $\tilde{\Sigma}$ of belief states, we first generated a set of initial distributions $\rho(x) = \phi(x; \mu_0, s^2)$ by fixing the variance to be $s^2 = 0.1$, and uniformly selecting the mean $\mu_0$ at random within the values of 17.5 and 22.  We then randomly sampled observations $y^x$ uniformly on $[16, 23.5]$, and chose an action $u\in\{0,1\}$ at random as well.  We continued this process for each $\sigma$ for $T$ time steps, and updated each $\sigma$ accordingly using \eqref{tempSigma}.  If a $\sigma$ function came too close to being everywhere zero, we reset that $\sigma$ to a new $\sigma_0$ and began the process again.  

We also used a mixture reduction process described in \cite{Zhang2010} to combine similar Gaussian components into a single new Gaussian based on the $L^2$ distance between the functions.  Once a new $\alpha$-function or $\sigma$ was generated from the previous $\alpha$-function or $\sigma$, we used the algorithm in \cite{Zhang2010} to reduce the total number of Gaussians to $20$. This helped reduce computation time, without overly sacrificing accuracy.  The number of components to keep can be easily changed, however, depending on the importance of trade-offs in speed versus accuracy.

\begin{figure*}[t!]
	\centering
	\subfloat[]{\label{probN5}\includegraphics[width=.5\textwidth]{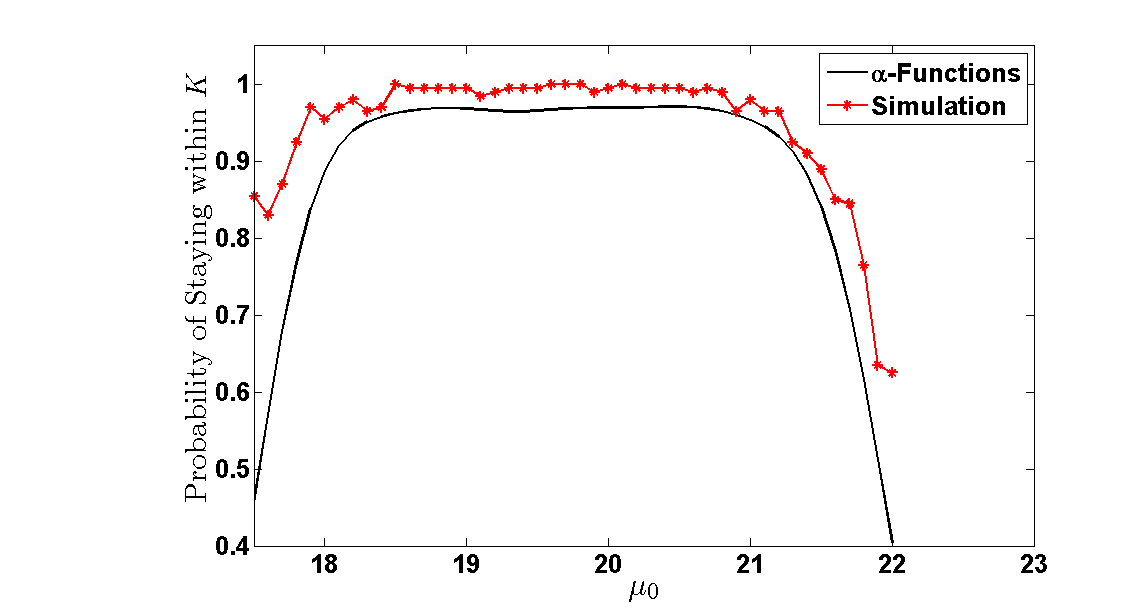}}
	\hfil
	\subfloat[]{\label{probN20}\includegraphics[width=.5\textwidth]{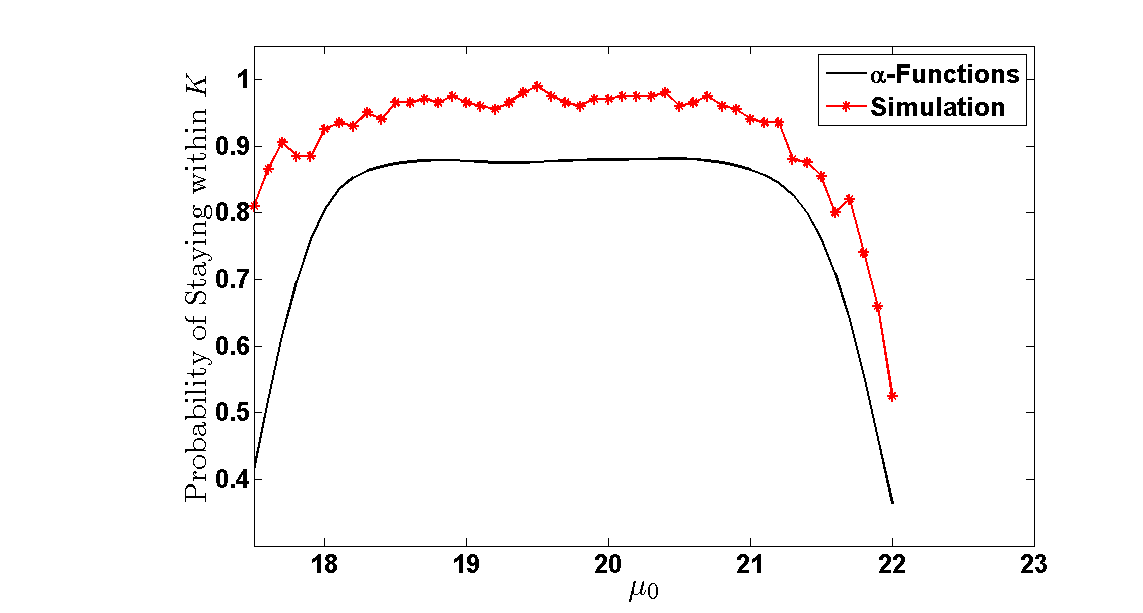}}

	\caption{Given $\sigma_0(x_0,q_0) = {\bf 1}_{\{0\}}(q_0)\phi(x_0;  \mu_0, 0.1)$, estimated probability of $x_k$ staying in $[17.5,22]$ for (a) $T=5$ time steps and (b) $T=20$ time steps according to $\alpha$-functions (in black) and according to simulation that uses the policy from the $\alpha$-functions (in red).  The $\alpha$-functions consistently underestimate the simulated reachability probability, assuring a minimum probability of safety, although the estimates from both have the same behavior, and are not too different except towards the boundaries of $K$.}
	\label{fig1}
\end{figure*}

Using a set $\tilde{\Sigma}$ of $40$ $\sigma$s, and running the backup operation $T$ times, we obtained an estimate of the probability of the temperature remaining within set $K$ for various $\mu_0$ values, and $q(0)=0$.  We also used the $\alpha$-functions calculated in the $T$th iteration as a stationary policy to estimate the average reachability probability for various $\mu_0$. To do so, we ran $200$ simulations of the temperature of the room over $T$  time steps for each $\mu_0$, using the stationary policy generated by the $\alpha$-functions to choose control actions.  The results of both the approximation to the probabilities via the value function estimate, as well as the probability estimated through simulation of the policy, are presented in Fig. \ref{fig1} for $T=5$  (\ref{probN5}) and $T=20$ (\ref{probN20}).  

\begin{figure}[!t]
 	\centering
  	\includegraphics[width=.5\textwidth]{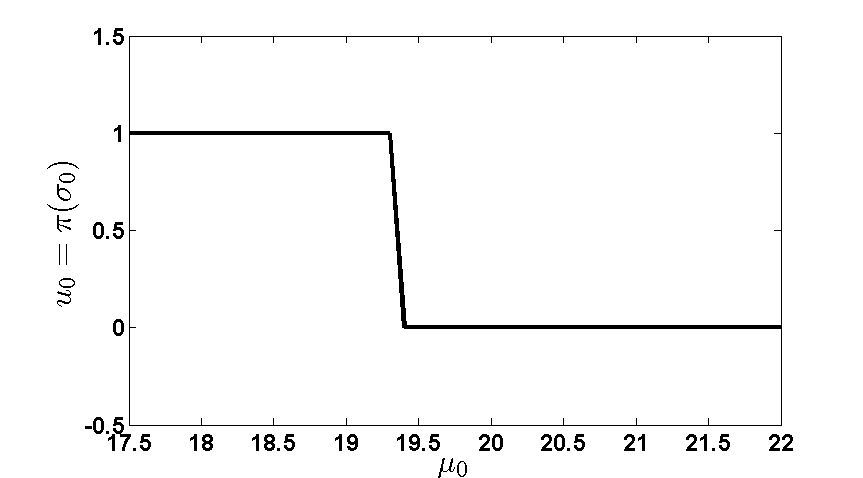}
	 \caption{Given $\sigma_0(x_0,q_0) = {\bf 1}_{\{0\}}(q_0)\phi(x_0;  \mu_0, 0.1)$, optimal choice of $u_0$ for varying $\mu_0$.  For $\mu_o$ less than or equal to 19.3, it is best to set $u_0=1$ (turn the heater on) and after 19.3, the heater should be left off.}
	 \label{fig2}
 \end{figure}

The value function estimate of the probability closely resembles the estimates from simulation, although near the edges of $K$ the discrepencies are larger.  This is likely due to the inaccuracies in the Gaussian sum approximation to the indicator function, which are much more noticeable at the boundaries.  The $\alpha$-functions also consistently produce lower probabilities than the simulated optimal policy, partly due to the inaccuracy of the indicator fuction approximation, but also because they are designed to produce a lower bound on the true value function. Fig. \ref{fig2} shows the optimal choice of $u_0$ according to the $\alpha$-functions for varying $\mu_0$ (i.e. for varying $\sigma_0$, since $\sigma_0\sim\mathcal{N}(\mu_0,s^2)$).  When the mean $\mu_0$ is less than or equal to 19.3, the heater should be turned on ($u_0=1$), and for larger values of $\mu_0$ the heater should remain off.  

Computation time to produce the $\alpha$-functions is intensive.  For $T=20$ time steps, generating the $\alpha$-functions took approximately eight hours to calculate on an Intel 3.40 GHz CORE i7-2600 CPU with 8 GB of RAM.  However, once the $\alpha$-functions have been calculated, using them to generate optimal control actions takes less than a second, including the time required to update the belief.  Thus, to estimate the probability of remaining in $K$ for a single sample trajectory, both when generating the $\sigma$ functions as Gaussian sums, and finding the optimal $\alpha$-function and associated control input, only takes a few seconds.

We found that the main bottleneck in computation was the discretized observations.  The PBVI algorithm requires evaluating the set of all observations at several different times in the backup process, and hence is not well suited to a large number of discrete observations.  Another issue in extending this mehod to higher dimensional systems is in approximating the indicator function as a sum of Gaussians, which is required both in the $\alpha$-function representation and in the belief update.  The growth in the number of Gaussians needed to adequately approximate a higher dimensional indicator function  slows down the overall computation time.  Therefore in order to apply this PBVI technique to higher dimensional systems, we will need to explore better representations of the indicator function (possibly using a particle filter to represent the beliefs, as in \cite{Porta06} or \cite{Thrun2000}) as well as efficient ways to allow for continuous observations, possibly by using the method described in \cite{Porta06}, which groups observations according to which $\alpha$-functions are optimal for those observations and creating a discrete representative for each group.

\section{Conclusion}\label{conc}

We have provided the first numerical results to the reachability problem with partially observable discrete time stochastic hybrid dynamics.  By showing that the value function is still piecewise-linear and convex in the case of discrete actions and observations, and that the representation of the $\alpha$-functions and belief states by linear combinations of Gaussians is preserved under the backup operator and belief update, we were able to extend PBVI techniques for continuous state systems to the reachability problem for PODTSHS.  We then demonstrated our method on a one dimensional temperature regulation problem with stochastic hybrid dynamics and a noisy discretized measurement of the continuous state.  Although the calculation of the $\alpha$-functions was slow, the policy they encode can be applied quickly online to optimize the system's probability of remaining within a safe region.  However, we hope to find more efficient techniques to overcome some of the current method's shortcomings.  Over larger state spaces, discretizing the observation space is not practical, and techniques that accomodate a continuous observation space should be explored.  The use of particle filters to represent the beliefs may also be beneficial, because of the inability of a small number of Gaussian components to adequately represent discontinuous functions (such as the indicator function).  Overall, we believe our method is well-suited to low dimensional systems, and with further investigation should be extendable to higher dimensional systems as well.

\bibliographystyle{IEEEtran}
\bibliography{POMDPbib}

\end{document}